\documentclass[12pt]{amsart}

\usepackage{amssymb}
\usepackage{bm}
\usepackage{graphicx}
\usepackage{psfrag}
\usepackage{hyperref}
\hypersetup{colorlinks=true, linkcolor=blue, citecolor=red, urlcolor=wine}
\usepackage{color}
\usepackage{url}
\usepackage{algpseudocode}
\usepackage{fancyhdr}
\usepackage{xy}
\input xy
\xyoption{all}
\usepackage{stmaryrd}
\usepackage{calrsfs}
\usepackage[dvipsnames]{xcolor}

\newtheorem*{maintheorem*}{Main Theorem}
\newtheorem{theorem}{Theorem}[section]
\newtheorem{prop}[theorem]{Proposition}

\newtheorem{conj}[theorem]{Conjecture}
\newtheorem{question}[theorem]{Question}
\newtheorem{lemma}[theorem]{Lemma}
\newtheorem{cor}[theorem]{Corollary}
\theoremstyle{definition}
\newtheorem{definition}[theorem]{Definition}
\newtheorem{remark}[theorem]{Remark}
\newtheorem{example}[theorem]{Example}
\numberwithin{equation}{section}

\newcommand{\ii}{\mathcal{A}}
\newcommand{\nn}{\mathbb{N}}
\newcommand{\pp}{\mathbb{P}}
\newcommand{\qq}{\mathbb{Q}}
\newcommand{\rr}{\mathbb{R}}
\newcommand{\uu}{\mathcal{U}}
\newcommand{\zz}{\mathbb{Z}}
\newcommand\pval{\mathsf{v}_p}
\newcommand{\gp}{\text{gp}}
\newcommand{\qf}{\text{qf}}
\newcommand{\supp}{\text{supp}}

\providecommand\ldb{\llbracket}
\providecommand\rdb{\rrbracket}

\keywords{atomic domain, BFD, bounded factorization domain, FFD, finite factorization domain, UFD, seminormal domain, isomorphism problem}


\begin{document}
	\mbox{}
	\title{On semigroup algebras with rational exponents}
	
	\author{Felix Gotti}
	\address{Department of Mathematics\\MIT\\Cambridge, MA 02139}
	\email{fgotti@mit.edu}
	
	\date{\today}
	
	\begin{abstract}
		In this paper, a semigroup algebra consisting of polynomial expressions with coefficients in a field $F$ and exponents in an additive submonoid $M$ of $\qq_{\ge 0}$ is called a Puiseux algebra and denoted by $F[M]$. Here we study the atomic structure of Puiseux algebras. To begin with, we answer the Isomorphism Problem for the class of Puiseux algebras, that is, we show that for a field $F$ if two Puiseux algebras $F[M_1]$ and $F[M_2]$ are isomorphic, then the monoids $M_1$ and $M_2$ must be isomorphic. Then we construct three classes of Puiseux algebras satisfying the following well-known atomic properties: the ACCP property, the bounded factorization property, and the finite factorization property. We show that there are bounded factorization Puiseux algebras with extremal systems of sets of lengths, which allows us to prove that Puiseux algebras cannot be determined up to isomorphism by their arithmetic of lengths. Finally, we give a full description of the seminormal closure, root closure, and complete integral closure of a Puiseux algebra, and use such description to provide a class of antimatter Puiseux algebras (i.e., Puiseux algebras containing no irreducibles).
	\end{abstract}

\maketitle

\section{Introduction} \label{sec:intro}

The study of group rings dates back to the mid-nineteen century. Most of the initial research in this area focused mainly on the structure of the groups of units of group rings and the Isomorphism Problem for group rings over a given ring of coefficients; see G. Highman~\cite{gH40} and S.~K. Sehgal~\cite{sS69}, respectively. However, most of the early ring-theoretical study of group rings was mainly carried out on a non-commutative setting. It was not until the seventies with the work of R. Gilmer, R.~Matsuda, and other authors that the study of commutative group rings as well as commutative semigroup rings started earning substantial attention (see \cite{rG77,GH79,rM79} and references therein). Much of the work done on commutative semigroup rings during this decade focused on the following abstract problem: given a commutative ring $R$ and a commutative semigroup $S$, establish conditions under which the semigroup ring of $S$ over $R$ satisfies certain algebraic property. Answering instances of this problem often requires a fair understanding of the algebraic properties of both $S$ and $R$. As the structure of commutative semigroups most of the time cannot be derived from that of abelian groups, a new research direction in commutative algebra had emerged.

Much of the work on commutative semigroup rings carried out in the seventies was compiled by Gilmer in his celebrated book Commutative Semigroup Rings~\cite{rG84}, which in turn has motivated a lot of research in the field. More recently, many authors, including P.~A. Grillet~\cite{pG95}, J. Gubeladze~\cite{jG98}, and H. Kim~\cite{hK01}, have investigated algebraic and factorization properties of semigroup rings. Semigroup algebras (i.e., semigroup rings with coefficients in a field), in particular, have permeated through various areas under active investigation, including algebraic combinatorics~\cite{BCMP98}, discrete geometry~\cite{BG09}, and functional analysis~\cite{mA04}.

Puiseux monoids (i.e., additive submonoids of nonnegative rationals) have a complex atomic structure~(see \cite{GGT21} and references therein). As a result, the semigroup algebras they determine, which we refer to as Puiseux algebras, have played important roles in commutative algebra. For instance, A. Grams in~\cite{aG74} localized a Puiseux algebra to construct an integral domain and disprove P.~M. Cohn's assumption that every atomic domain satisfies the ACCP. In addition, J. Coykendall and the author recently appealed to Puiseux algebras~\cite{CG19} to partially answer a question on atomicity stated by Gilmer back in 1984 \cite[page~189]{rG84}. There are further appearances of Puiseux algebras in recent literature (for instance, in~\cite{ACHZ07} and~\cite{GK20}); however, no systematic study of their atomic structure seems to have been carried out so far. Although this paper offers by no means a systematic study, it aims to provide further insight on the algebraic and atomic structure of Puiseux algebras. Here we address two algebraic problems: the Isomorphism Problem and the computation of the seminormal, root, and complete integral closures for Puiseux algebras. Then we use both results to construct various infinite classes of Puiseux algebras with distinct atomic properties.

The first problem we shall address here is the Isomorphism Problem for Puiseux algebras. The Isomorphism Problem for a field $F$ and a class of monoids $\mathcal{C}$ is the question of whether two monoids in $\mathcal{C}$ are isomorphic provided they have isomorphic semigroup algebra over $F$. Versions of the Isomorphism Problem on classes of finitely generated monoids have been investigated before; see, for instance,~\cite{pG95} and~\cite{jG98}. However, the Isomorphism Problem on classes of non-finitely generated monoids seems to be rather unexplored. In Section~\ref{sec:the Isomorphism Problem}, we give a positive answer to the Isomorphism Problem for Puiseux algebras.

Let $R$ be an integral domain. The domain $R$ is said to satisfy the ascending chain condition on principal ideals (ACCP) if every ascending chain of principal ideals of $R$ becomes stationary. In addition, $R$ is called a bounded factorization domain (BFD) if for every nonzero nonunit $x \in R$ there exists $N \in \nn$ such that $x = a_1 \cdots a_n$ for irreducibles $a_1, \dots, a_n \in R$ implies that $n \le N$. Moreover, $R$ is called a finite factorization domain (FFD) if every nonzero element of $R$ has only finitely many non-associate divisors. The notions of BFDs and FFDs were introduced in~\cite{AAZ90} by D. D. Anderson, D. F. Anderson, and M. Zafrullah, where the authors introduced and studied a diagram of implications of atomic classes of integral domains containing the chain
\[
	\textbf{UFD} \Longrightarrow \ \textbf{FFD} \Longrightarrow \ \textbf{BFD} \Longrightarrow \ \textbf{ACCP}.
\]
These three implications are not reversible in general, and examples of atomic domains witnessing this observation are provided in~\cite{AAZ90}. To illustrate the complexity of the class of Puiseux algebras, for each of the three implications, we construct in Section~\ref{sec:atomic PA} an infinite class of Puiseux algebras witnessing the failure of its converse. Our positive answer to the Isomorphism Theorem is key to guarantee that the classes we construct are infinite up to isomorphism.

For each nonzero element $x$ of an atomic integral domain $R$, the set $\mathsf{L}(x)$ consists of all possible lengths of factorizations of $x$, and $\mathcal{L}(R) := \{\mathsf{L}(x) \mid x \in R \setminus \{0\}\}$ is called the system of sets of lengths of $R$; the system of sets of lengths of an atomic monoid is defined similarly. The search for bounded factorization domains and monoids having extremal systems of sets of lengths has earned significant interest since F. Kainrath~\cite{fK99} proved that the system of sets of lengths of a Krull domain/monoid with infinite class group (and primes in each divisor class) is as large as it can possible be (see, for instance, \cite{FNR19,fG19a}). In Section~\ref{sec:atomic PA} we show that there are infinitely many non-isomorphic Puiseux algebras having extremal systems of sets of lengths. This will allow us to answer negatively the Characterization Problem for Puiseux algebras: for any fixed field $F$, Puiseux algebras over $F$ cannot be determined up to isomorphism by their systems of sets of lengths.

Following Coykendall et al.~\cite{CDM99}, we say that an integral domain is antimatter if it contains no irreducibles. Classes of antimatter Puiseux algebras have been constructed in~\cite{ACHZ07}. Following A.~Grams and H.~Warner~\cite{GW75}, we say that an integral domain is irreducible-divisor-finite (or IDF) if each element is divisible by only finitely many irreducibles up to associates; IDF monoids are defined similarly. Every antimatter domain/monoid is clearly IDF, and an integral domain (resp., a monoid) satisfies the finite factorization property if and only if it is atomic and IDF by \cite[Theorem~5.1]{AAZ90} (resp., \cite[Theorem~2]{fHK92}). In Section~\ref{sec:antimatter PA}, we show that a Puiseux algebra $F[M]$ may not be IDF even when $M$ is IDF. In Section~\ref{sec:antimatter PA}, we also prove that the seminormal closure, root closure, and complete integral closure of a Puiseux algebra all coincide, and we provide a full description in terms of its monoid of exponents. We use our description to determine whether the seminormal closure of a Puiseux algebra is atomic or antimatter. In sharp contrast to the fact that the Puiseux algebras we obtain by taking seminormal closures are seminormal, we conclude the paper providing a non-seminormal class of antimatter Puiseux algebras.

\bigskip
\section{Background} 
\label{sec:background}

Throughout this paper, we let $\mathbb{N}$ and $\pp$ denote the set of positive integers and the set of primes, respectively. In addition, we set $\nn_0 := \nn \cup \{0\}$. If $j,k \in \zz$, then we let $\ldb j,k \rdb$ denote the discrete interval from $j$ to $k$, i.e., $\ldb j,k \rdb := \{z \in \zz \mid j \le z \le k\}$. For $X \subseteq \rr$ and $r \in \rr$, we define $X_{\ge r} := \{x \in X \mid x \ge r\}$ and, in a similar manner, we define $X_{< r}$ and $X_{> r}$. If $q \in \qq_{> 0}$, then we denote the unique $n,d \in \nn$ such that $q = n/d$ and $\gcd(n,d)=1$ by $\mathsf{n}(q)$ and $\mathsf{d}(q)$, respectively. Finally, for $S \subseteq \qq_{> 0}$, we set
\[
	\mathsf{n}(S) := \{\mathsf{n}(s) \mid s \in S\} \quad \text{ and } \quad \mathsf{d}(S) := \{\mathsf{d}(s) \mid s \in S\}.
\]

\medskip
\subsection{Commutative Monoids}

Within the scope of our exposition, each monoid is tacitly assumed to be cancellative and commutative. Unless we specify otherwise, monoids here are written additively. Let~$M$ be a monoid. We set $M^\bullet := M \setminus \{0\}$, and we let $\uu(M)$ denote the set of units (i.e., invertible elements) of $M$. If $\uu(M) = \{0\}$, then $M$ is called \emph{reduced}. The quotient monoid $M/\uu(M)$, denoted by~$M_{\text{red}}$, is clearly reduced. If $x,y \in M$, then $y$ \emph{divides} $x$ \emph{in} $M$, in symbols, $y \mid_M x$, if there exists $z \in M$ such that $x = y + z$.
\smallskip

The \emph{difference group} of $M$, denoted here by $\gp(M)$, is the abelian group (unique up to isomorphism) satisfying that any abelian group containing a homomorphic image of~$M$ also contains a homomorphic image of $\gp(M)$. We say that $y \in M^\bullet$ has \emph{type zero} provided that there is a largest $n \in \nn$ such that the equation $nx = y$ is solvable in $\gp(M)$. In addition, we say that $M$ has \emph{type zero} if every element of $M^\bullet$ has type zero. The monoids
\begin{itemize}
	\item \hspace{-6pt} $M' := \big\{ x \in \gp(M) \mid \text{ there exists } N \in \nn \text{ such that } nx \in M \ \text{for all} \ n \ge N \big\}$,
	
	\item $\bar{M} := \big\{ x \in \gp(M) \mid \ nx \in M \ \text{for some} \ n \in \nn \big\}$, and
	
	\item $\widehat{M} := \big\{ x \in \gp(M) \mid \text{ there exists } c \in M \text{ such that } c + nx \in M \text{ for all } n \in \nn \big\}$
\end{itemize}
are called the \emph{seminormal closure, root closure}, and \emph{complete integral closure} of $M$, respectively. In addition, the following chain of inclusions holds:
\begin{equation} \label{eq:inclusion chain of closures}
	M \subseteq M' \subseteq \bar{M} \subseteq \widehat{M} \subseteq \gp(M).
\end{equation}
The only inclusion that is not obvious in \eqref{eq:inclusion chain of closures} is $\bar{M} \subseteq \widehat{M}$. To verify that this inclusion holds, take $g = x-y \in \bar{M}$ for some $x,y \in M$, and then take $n_0 \in \nn$ with $n_0 g \in M$ and set $c = n_0 y$. Now for every $n \in \nn$, we can write $n = m n_0 + r$ for some $m \in \nn_0$ and $r \in \ldb 0, n_0-1 \rdb$ to see that
\begin{equation} \label{eq:root closure inside complete integral closure}
	c + ng = n_0 y + (m n_0 + r)(x-y) = (n_0 - r)y + rx + m(n_0 g) \in M,
\end{equation}
whence $g \in \widehat{M}$. As a consequence, $ \bar{M} \subseteq \widehat{M}$. The monoid~$M$ is called \emph{seminormal} (resp., \emph{root closed} or \emph{completely integrally closed}) provided that the equality $M' = M$ (resp., $\bar{M} = M$ or $\widehat{M} = M$) holds.
\smallskip

An element $a \in M \setminus \uu(M)$ is called an \emph{atom} if for all $x,y \in M$ with $a = x+y$ either $x \in \uu(M)$ or $y \in \uu(M)$. The set of atoms of~$M$ is denoted by $\mathcal{A}(M)$. Atomicity and antimatterness play a fundamental role in this paper.

\begin{definition}
	Let $M$ be a monoid.
	\begin{enumerate}
		\item If each nonunit of $M$ can be written as a sum of atoms, then $M$ is \emph{atomic}.
		\smallskip
		
		\item If $M$ contains no atoms, then $M$ is \emph{antimatter}.
	\end{enumerate}
\end{definition} 

\noindent A subset $I$ of $M$ is an \emph{ideal} of $M$ provided that $I + M \subseteq I$. The ideal $I$ is \emph{principal} if $I = x + M$ for some $x \in M$. The monoid $M$ satisfies the \emph{ascending chain condition on principal ideals} (or \emph{ACCP}) if each increasing sequence of principal ideals of~$M$ eventually stabilizes. It is well known that each monoid satisfying the ACCP is atomic \cite[Proposition~1.1.4]{GH06b}. For $S \subseteq M$, we let $\langle S \rangle$ denote the submonoid of $M$ generated by~$S$ (i.e., the smallest submonoid of $M$ containing~$S$). The monoid ~$M$ is called \emph{finitely generated} if it can be generated by a finite set, while $M$ is called \emph{cyclic} if it can be generated by a singleton. Each finitely generated monoid satisfies the ACCP \cite[Proposition~2.7.8]{GH06b}.

In this paper we study monoid algebras whose exponents lie in $\qq_{\ge 0}$. Additive submonoids of $\qq_{\ge 0}$ are known as \emph{Puiseux monoids}. Additive submonoids of $\qq$ account, up to isomorphism, for all rank-one torsion-free monoids~\cite[Section~24]{lF70}, and every additive submonoid of $\qq$ that is not a group is isomorphic to a Puiseux monoid~\cite[Theorem~2.9]{rG84}. The atomic spectrum of the class of Puiseux monoids is broad, whose members ranging from antimatter monoids (e.g., $\langle 1/2^n \mid n \in \nn \rangle$) to non-finitely generated atomic monoids (e.g., $\langle 1/p \mid p \in \pp \rangle$). The atomic structure of Puiseux monoids has been systematically studied during the last three years, and the most relevant achieved results can be found in the survey~\cite{CGG20a}. Puiseux monoids have also been studied in connection with factorizations of matrices~\cite{BG20} and commutative rings~\cite{CG19}.

\medskip
\subsection{Factorizations}

Recall that a multiplicative monoid $F$ is \emph{the} free commutative monoid on $P \subseteq F$ provided that every element $x \in F$ can be written uniquely in the form $x = \prod_{p \in P} p^{\pval(x)}$, where $\pval(x) \in \nn_0$ and $\pval(x) > 0$ only for finitely many $p \in P$. For each set $P$, there exists a unique (up to canonical isomorphism) monoid $F$ that is free commutative on~$P$. We let $\mathsf{Z}(M)$ denote the (multiplicative) free commutative monoid on $\mathcal{A}(M_{\text{red}})$. The elements of $\mathsf{Z}(M)$ are called \emph{factorizations}. If $z = a_1 \cdots a_\ell \in \mathsf{Z}(M)$ for some $\ell \in \nn_0$ and $a_1, \dots, a_\ell \in \mathcal{A}(M_{\text{red}})$, then $\ell$ is called the \emph{length} of $z$ and is denoted by $|z|$. Because $\mathsf{Z}(M)$ is free, there exists a unique monoid homomorphism $\pi \colon \mathsf{Z}(M) \to M_{\text{red}}$ satisfying that $\pi(a) = a$ for all $a \in \mathcal{A}(M_{\text{red}})$. For each $x \in M$, we set
\[
	\mathsf{Z}(x) := \mathsf{Z}_M(x) := \pi^{-1}(x + \uu(M)) \subseteq \mathsf{Z}(M).
\]
Clearly, $M$ is an atomic monoid if and only if $\mathsf{Z}(x)$ is nonempty for all $x \in M$. The monoid $M$ is called a \emph{unique factorization monoid} (or a \emph{UFM}) if $|\mathsf{Z}(x)| = 1$ for all $x \in M$. On the other hand, $M$ is called a \emph{finite factorization monoid} (or an \emph{FFM}) if $1 \le |\mathsf{Z}(x)| < \infty$ for all $x \in M$. Clearly, each UFM is an FFM. For each $x \in M$, we set
\[
	\mathsf{L}(x) := \mathsf{L}_M(x) := \{|z| \ \mid z \in \mathsf{Z}(x)\}.
\]
The set $\mathcal{L}(M) := \{\mathsf{L}(x) \mid x \in M\}$ is called the \emph{system of sets of lengths} of $M$. The monoid $M$ is called a \emph{bounded factorization monoid} (or a \emph{BFM}) if $1 \le |\mathsf{L}(b)| < \infty$ for all $b \in M$. It is clear that each FFM is a BFM. In addition, it is well known and easy to prove that each BFM satisfies the ACCP (see the Diagram~\eqref{diag:AAZ chain}). Sets of lengths play an important role in the study of the arithmetic of atomic monoids (see~\cite{aG16} for more details).

\medskip
\subsection{Integral Domains and Semigroup Rings} 
\label{sec:semigroup algebras of PM}

Let $R$ be an integral domain, and let $R^* := R \setminus \{0\}$ denote the multiplicative monoid of $R$. As usual, $R^\times$ denote the group of units of $R$ (clearly, $R^\times = \uu(R^*)$). We say that~$R$ is \emph{atomic} (resp., \emph{antimatter}, a \emph{BFD}, an \emph{FFD}) if $R^*$ is atomic (resp., antimatter, a BFM, an FFM). For several examples of BFDs and FFDs, see the recent survey~\cite{AG20}. Also, we let $\mathcal{A}(R)$ and $\mathsf{Z}(R)$ denote $\mathcal{A}(R^*)$ and $\mathsf{Z}(R^*)$, respectively, and for a nonzero nonunit $r \in R$, we let $\mathsf{Z}(r) := \mathsf{Z}_R(r)$ and $\mathsf{L}(r) := \mathsf{L}_R(r)$ denote $\mathsf{Z}_{R^*}(r)$ and $\mathsf{L}_{R^*}(r)$, respectively. The \emph{system of sets of lengths} of~$R$ is defined to be $\mathcal{L}(R) := \mathcal{L}(R^\ast)$

For a ring $R$ and a semigroup $S$, consider the set $R[X;S]$ comprising all functions $f \colon S \to R$ satisfying that $\{s \in S \mid f(s) \neq 0 \}$ is finite. We shall conveniently represent an element $f \in R[X;S]$ by $f = \sum_{i=1}^n f(s_i)X^{s_i}$, where $s_1, \dots, s_n$ are precisely those $s \in S$ satisfying that $f(s) \neq 0$. With addition and multiplication defined as for polynomials, $R[X;S]$ is a ring, which is called the \emph{semigroup ring of} $S$ \emph{over} $R$. Following Gilmer~\cite{rG84}, we shall write $R[S]$ instead of $R[X;S]$. As we are mainly concerned with semigroup rings of Puiseux monoids over a given field, the following terminology seems appropriate.

\begin{definition}
	Let $F$ be a field, and let $M$ be a Puiseux monoid. Then we call $F[M]$ a \emph{Puiseux algebra}.
\end{definition}

Let $F[M]$ be a Puiseux algebra. It follows from \cite[Theorem~8.1 and Theorem~11.1]{rG84} that $F[M]$ is an integral domain satisfying that $F[M]^\times = F$. We say that the element $f = \alpha_1X^{q_1} + \dots + \alpha_k X^{q_k} \in F[M]^*$ is represented in \emph{canonical form} if $\alpha_i \neq 0$ for each $i \in \ldb 1,k \rdb$ and $q_1 > \dots > q_k$. Observe that each element of $F[M]^*$ has a unique representation in canonical form. In this case, $\supp(f) := \{q_1, \dots, q_k\}$ and $\deg(f) := q_1$ are called the \emph{support} and the \emph{degree} of $f$, respectively. As for polynomials, the degree identity $\deg fg = \deg f + \deg g$ holds for all $f,g \in F[M]^*$.

\bigskip
\section{The Isomorphism Problem}
\label{sec:the Isomorphism Problem}

For each integral domain $R$, a monoid isomorphism $M_1 \to M_2$ always induces an $R$-algebra isomorphism $R[M_1] \to R[M_2]$. More generally, we have the following proposition, which is an immediate consequence of~\cite[Theorem~7.2(2)]{rG84}.

\begin{prop} \label{prop:induced homomorphism of semigroup algebras}
	Let $R$ be an integral domain, and let $\varphi \colon M_1 \to M_2$ be a monoid homomorphism. If $\varphi$ is injective (resp., surjective), then the $R$-algebra homomorphism $\bar \varphi \colon R[M_1] \to R[M_2]$ determined by $X^s \mapsto X^{\varphi(s)}$ is also injective (resp., surjective).
\end{prop}

In the context of monoid algebras, the Isomorphism Problem refers to the veracity of the reverse implication of Proposition~\ref{prop:induced homomorphism of semigroup algebras}.
\medskip

\noindent \textbf{Isomorphism Problem.} Let $F$ be a field, and let $\mathcal{C}$ be a nonempty class of monoids. For $M_1, M_2 \in \mathcal{C}$, does $F[M_1] \cong F[M_2]$ as $F$-algebras guarantee that $M_1 \cong M_2$?
\medskip

A brief early survey on this problem is offered by Gilmer in~\cite[Section~25]{rG84}. Although the cases when $\mathcal{C}$ is a class consisting of groups have been studied since the 1960s~\cite{sB67}, it was not until the 1980s that some attention was given to the more general case of monoids. In 1982, A. S. Demushkin proposed in~\cite{aD82} a positive answer to the Isomorphism Problem on the class of finitely generated normal monoids without nontrivial units (see~\cite[Section~11]{rG84} for definitions). However, his proof involved various invalid arguments. In 1998, Gubeladze provided a final positive answer to the Isomorphism Problem on the class of finitely generated torsion-free monoids~\cite{jG98}. We proceed to offer a positive answer to the Isomorphism Problem on the class of Puiseux monoids.

\begin{theorem} \label{thm:the isomorphism problem}
	Let $F$ be a field, and let $M_1$ and $M_2$ be Puiseux monoids. Then $M_1 \cong M_2$ as monoids if and only if $F[M_1] \cong F[M_2]$ as $F$-algebras.
\end{theorem}

\begin{proof}
	The direct implication is an immediate consequence of Proposition~\ref{prop:induced homomorphism of semigroup algebras}. To prove the reverse implication, suppose that $\bar \varphi \colon F[M_1] \to F[M_2]$ is an $F$-algebra isomorphism. By virtue of Proposition~\ref{prop:induced homomorphism of semigroup algebras} one can replace $M_1$ by an isomorphic copy $rM_1$ (for a suitable $r \in \qq_{> 0}$), and therefore, assume that both conditions $F[M_1] \cong F[M_2]$ and $X \in F[M_1]$ hold. Now define $\varphi \colon M_1 \to M_2$ by $\varphi(q) = \deg \bar\varphi(X^q)$, and notice that for each $q \in M_1^\bullet$,
	\[
		\varphi(q) = \deg \bar{\varphi}(X^q) = \frac{1}{\mathsf{d}(q)} \deg \bar{\varphi}(X^q)^{\mathsf{d}(q)} = \frac{1}{\mathsf{d}(q)} \deg \bar{\varphi}(X)^{\mathsf{n}(q)} = q \deg \bar{\varphi}(X).
	\]
	Hence after setting $q_0 := \deg \bar\varphi(X) \in M_2 \subseteq \qq_{\ge 0}$, we see that $\varphi$ is the monoid homomorphism consisting in multiplying by $q_0$. Since $\bar \varphi$ is an $F$-algebra isomorphism, $q_0 > 0$. Therefore $\varphi$ is not only injective, but also strictly increasing. Finally, we show that $\varphi$ is surjective. Since $\varphi$ is strictly increasing, for each element $f := \sum_{k=1}^n c_k X^{q_k} \in F[M_1]^\ast$  represented in canonical form, one obtains that $\deg \bar\varphi(X^{q_n}) > \deg \bar\varphi(X^{q_k})$ for every $k \in \ldb1, n-1 \rdb$, and as a consequence,
	\begin{equation} \label{eq:degrees}
		\deg \bar\varphi(f) = \deg \, \sum_{k=1}^n c_k \bar\varphi \big(X^{q_k} \big) = \deg \bar\varphi\big( X^{q_n}\big) = \varphi(q_n) = q_0 \deg f.
	\end{equation}
	Clearly, $M_1^\bullet = \{\deg f \mid f \in F[M_1] \setminus F \}$ and $M_2^\bullet = \{\deg g \mid g \in F[M_2] \setminus F \}$. As $\bar\varphi$ is a surjective function satisfying that $\bar \varphi(F) = F$, it follows from~\eqref{eq:degrees} that
	\[
		M_2^\bullet = \{\deg \bar\varphi (f) \mid f \in F[M_1] \setminus F\} = q_0 \{\deg f \mid f \in F[M_1] \setminus F\} = q_0 M_1^\bullet.
	\]
	As $q_0M_1 = M_2$, the homomorphism $\varphi$ is surjective. Hence $M_1 \cong M_2$ as monoids.
\end{proof}

\bigskip
\section{Classes of Atomic Puiseux Algebras}
\label{sec:atomic PA}

The chain~\eqref{diag:AAZ chain} of refined classes of atomic domains was introduced by Anderson, Anderson, and Zafrullah in~\cite{AAZ90}. Since then this chain has received a significant amount of consideration in the literature of both commutative algebra and semigroup theory.
\begin{equation} \label{diag:AAZ chain}
	\textbf{UFD} \Longrightarrow \ \textbf{FFD} \Longrightarrow \ \textbf{BFD} \Longrightarrow \ \textbf{ACCP} \Longrightarrow \ \textbf{atomic domain}
\end{equation}
As illustrated in~\cite{AAZ90}, none of the implications in~\eqref{diag:AAZ chain} is, in general, reversible. This section is devoted to study the potential failure of each of the reverse implications in~\eqref{diag:AAZ chain} when we restrict to the class of Puiseux algebras. We will construct atomic Puiseux algebras witnessing such failure for the three leftmost implications, illustrating, as a byproduct, the diversity and complexity of the atomic structure of Puiseux algebras. We still do not know whether there exists an atomic Puiseux algebra failing to satisfy the ACCP. However, we suspect that this is the case, and we propose a potential witness at the end of this section.

\begin{remark}
	The original full diagram containing the chain~\eqref{diag:AAZ chain} also involves the class of half-factorial domains. An integral domain $R$ is called \emph{half-factorial} (or an \emph{HFD}) provided that $|\mathsf{L}(x)| = 1$ for every $x \in R^*$. We shall not explicitly consider half-factoriality here because a Puiseux algebra is an HFD if and only if it is a UFD \cite[Theorem~4.4]{fG20}.
\end{remark}

If a Puiseux algebra $F[M]$ is atomic, then so is the monoid $M$~\cite[Proposition~1.4]{hK01}. The converse statement was posed by Gilmer in~\cite[page~189]{rG84} and has been answered negatively by Coykendall and the author~\cite[Theorem~5.4]{CG19}. In this section, we identify various infinite classes of atomic Puiseux monoids whose corresponding Puiseux algebras are also atomic but play different roles in the chain of atomic classes~\eqref{diag:AAZ chain}.

\medskip
\subsection{The Ascending Chain Condition on Principal Ideals} 

If $R$ is an integral domain and~$M$ is a torsion-free monoid such that either $M$ is reduced or has type zero, then $R[M]$ satisfies the ACCP if and only if both $R$ and $M$ satisfy the ACCP (this was first noted by R. Gilmer and T. Parker in~\cite[Section 7]{GP74}). Using this result, we can construct classes of Puiseux algebras that satisfy the ACCP but are not BFDs. Part~(2) of Proposition~\ref{prop:PAlgebra ACCP transfer} is a slight generalization of \cite[Example~2.1]{AAZ90}, and its proof follows the same argument; we have included it here for the sake of completeness.

\begin{prop} \label{prop:PAlgebra ACCP transfer} \label{prop:PA that is ACCP but not BF}
	For a field $F$ the following statements hold.
	\begin{enumerate}
		\item If $M$ is a Puiseux monoid satisfying the ACCP, then the Puiseux algebra $F[M]$ also satisfies the ACCP.
		\smallskip
		
		\item  If $P \subseteq \pp$ is nonempty and $M_P = \langle 1/p \mid p \in P \rangle$, then $F[M_P]$ satisfies the ACCP, but it is a BFD if and only if $|P| < \infty$. In addition, for any two distinct  infinite subsets~$P$ and~$Q$ of $\pp$, it follows that $F[M_P] \cong F[M_Q]$ if and only if $P = Q$.
	\end{enumerate}
\end{prop}

\begin{proof}
	(1) This can be proved easily, but it also follows from the previously mentioned Gilmer and Parker's observation because every field trivially satisfies the ACCP and every Puiseux monoid is reduced.
	\smallskip
	
	(2) Fix a nonempty $P \subseteq \pp$, and let $M_P$ be as in part~(2). If $|P| < \infty$, then $M_P$ is an FFM by \cite[Proposition~2.7.8]{GH06b}, and so it satisfies the ACCP. We assume, therefore, that $|P| = \infty$. It is easy to check that $M_P$ is atomic with $\mathcal{A}(M_P) = \{1/p \mid p \in P\}$. Let $(p_n)_{n \in \nn}$ be a strictly increasing sequence with underlying set $P$. One can readily check that for each $q \in M$ there is a unique $N(q) \in \nn_0$ and a unique sequence of nonnegative integers $(c_n(q))_{n \in \nn}$ such that $q = N(q) + \sum_{n \in \nn} c_n(q) \frac{1}{p_n}$, where $c_n(q) \in \ldb 0, p_n - 1 \rdb$ and $c_n(q) = 0$ for all but finitely many $n \in \nn$. Set $s(q) := \sum_{n \in \nn} c_n(q)$. Clearly, if $q' \mid_M q$ for some $q' \in M$, then $N(q') \le N(q)$. In addition, observe that if $q'$ is a proper divisor of $q$ in $M$, then $N(q') = N(q)$ implies that $s(q') < s(q)$. As a consequence of these two observations, one deduces that each sequence $(q_n)_{n \in \nn}$ in $M$ satisfying that $q_{n+1} \mid_M q_n$ for every $n \in \nn$ must stabilize. Hence $M$ satisfies the ACCP.
	
	 We verify now that $F[M_P]$ is not a BFD when $|P| = \infty$. First notice that $a \in \mathcal{A}(M_P)$ if and only if $X^a \in \ii(F[M_P])$. Hence $X^{1/p} \in \ii(F[M_P])$ for every $p \in P$. Because $X = \big( X^{1/p} \big)^p$, it follows that $\mathsf{L}_{F[M_P]}(X) = P$. Since $|P| = \infty$, the Puiseux algebra $F[M_P]$ is not a BFD.
	 
	To argue the last statement, suppose that $P$ and $Q$ are infinite subsets of $\pp$ such that $F[M_P] \cong F[M_Q]$. It follows from Theorem~\ref{thm:the isomorphism problem} that $M_P \cong M_Q$, and then it follows from \cite[Proposition~3.2]{fG18} that $M_Q = rM_P$ for some $r \in \qq_{> 0}$. If $p \in P$, then $r/p \in \mathcal{A}(M_Q)$ and so $r = p/q$ for some $q \in Q$. Similarly, if $p_1 \in P \setminus \{p\}$, then $r = p_1/q_1$ for some $q_1 \in Q$. The equality $p q_1 = p_1 q$ now implies that $p = q \in Q$. Thus, $P \subseteq Q$. In a similar manner, one can argue that $Q \subseteq P$. The reverse implication is obvious.
\end{proof}

\medskip
\subsection{The Bounded Factorization Property}

For an integral domain $R$ and a torsion-free monoid $M$ that is reduced or has type zero, $R[M]$ is a BFD if and only if $R$ is a BFD and $M$ is a BFM. This was proved by H. Kim \cite[Theorem~3.15]{hK98} for monoids of type zero and by D. D. Anderson and J.~Juett \cite[Theorem~13]{AJ15} for reduced monoids. We will use the later result in Proposition~\ref{prop:BF for PAlgebras} to construct an infinite class of Puiseux algebras that are BFDs but not FFDs.

For every $p \in \pp$, the nonnegative cone of the localization of the ring $\zz$ at its multiplicative subset $\{p^n \mid n \in \nn_0\}$ is the antimatter Puiseux monoid $\langle \frac{1}{p^n} \mid n \in \nn_0 \rangle$. More generally, we can consider multiplicative subsets of $\zz$ generated by several primes, and take positive rays of their localization rings to obtain Puiseux monoids that are indeed atomic. The following example illustrates this observation.

\begin{example} \label{ex:Puiseux monoid by localizing Z}
	Let $P$ be the multiplicative subset of $\zz$ generated by the primes $2$ and~$3$, namely, $P = \{2^i 3^j \mid i,j \in \nn_0\}$. Then we see that
	\[
		\zz P^{-1}_{\ge 1} = \Big\{ \frac{n}{2^i 3^j} \ \big{|} \ n,i,j \in \nn_0 \Big\}_{\ge 1},
	\]
	and so $\{0\} \cup \zz P^{-1}_{\ge 1} = \{0\} \cup M_{\ge 1}$, where $M = \big\langle \frac{1}{2^i 3^j} \mid i,j \in \nn_0 \big\rangle$. Observe that although~$M$ is antimatter, the Puiseux monoid $\{0\} \cup \zz P^{-1}_{\ge 1}$ is atomic; indeed, it is a BFM because~$0$ is not a limit point of $\zz P^{-1}_{\ge 1}$. On the other hand, it is easy to check that $3$ has infinitely many divisors in $\{0\} \cup \zz P^{-1}_{\ge 1}$, and so the later is not an FFM.
\end{example}

A generalized version of the monoid in Example~\ref{ex:Puiseux monoid by localizing Z} can be used to construct a class of Puiseux algebras that are BFDs but not FFDs.

\begin{prop} \label{prop:BF for PAlgebras}
	Let $F$ be a field, $\zz P^{-1}$ the localization of the ring $\zz$ at a multiplicative subset $P$ generated by primes, and $M_P = \{0\} \cup \zz P^{-1}_{\ge 1}$. Then the following statements hold.
	\begin{enumerate}
		\item If a Puiseux monoid $M$ is a BFM, then the Puiseux algebra $F[M]$ is a BFD.
		\smallskip
		
		\item The Puiseux algebra $F[M_P]$ is a BFD, but it is an FFD if and only if $P = \{1\}$. In addition, $F[M_P] \cong F[M_Q]$ if and only if  $P = Q$.
	\end{enumerate}
\end{prop}

\begin{proof}
	(1) As $M$ is reduced, it is a direct consequence of part~(3) of \cite[Theorem 13]{AJ15}.
	\smallskip
	
	(2) Since $0$ is not a limit point of $M_P^\bullet$, it follows from \cite[Proposition~4.5]{fG19} that $M_P$ is a BFM. Hence $F[M_P]$ is a BFD by part~(1).
	
	If $P = \{1\}$, then $F[M_P] = F[\nn_0] = F[X]$, which is clearly an FFD. Now suppose that $P \neq \{1\}$, and let us verify that $F[M_P]$ is not an FFD. One can readily check that $\mathcal{A}(M_P) = [1, 2) \cap \zz P^{-1}$ and, as a consequence, the set of irreducibles in $F[M_P]$ dividing $X^3$ is $A = \{X^q \mid q \in [1,2) \cap \zz P^{-1} \}$. Let $p$ be one of the primes generating $P$. Because the irreducibles in $A$ are pairwise non-associate, the equalities $X^3 = X^{1 + 1/p^n} X^{2 - 1/p^n}$ (for every $n \in \nn$) yield infinitely many factorizations of $X^3$ in $F[M_P]$. Hence $F[M_P]$ is not an FFD.
	
	For the direct implication of the last statement, suppose that $F[M_P] \cong F[M_Q]$, where~$Q$ is also a multiplicative set of $\zz$ generated by primes. Using Theorem~\ref{thm:the isomorphism problem} we obtain that $M_P \cong M_Q$ and, therefore, \cite[Proposition~3.2]{fG18} guarantees that $M_Q = rM_P$ for some $r \in \qq_{> 0}$. As $\inf M_P^\bullet = \inf M_Q^\bullet = 1$, we see that $r=1$. Hence $M_P = M_Q$, and so $P=Q$. The reverse implication of the last statement is obvious.
\end{proof}
\smallskip

An integral domain (resp., a monoid) is said to have \emph{full system of sets of lengths} if it is a BFD (resp., BFM) and each subset of $\nn_{\ge 2}$ is a set of lengths of some element (i.e., its system of sets of lengths is as large as it can be). In the next proposition, we show that there are infinitely many non-isomorphic Puiseux algebras having full systems of sets of lengths. First, we need the following lemma.

\begin{lemma} \label{lem:PMs with full systems of sets of lengths}
	There are infinitely many non-isomorphic Puiseux monoids having full systems of sets of lengths.
\end{lemma}

\begin{proof}
	For each infinite set of primes $P$, one can mimic the proof of \cite[Theorem~3.6]{fG19a} to construct a non-finitely generated Puiseux monoid $M$ having full system of sets of lengths such that $\mathsf{d}(M^\bullet)$ is contained in the free (multiplicative) monoid with base $P$. Now let $(P_n)_{n \in \nn}$ be a sequence of pairwise disjoint subsets of $\pp$ such that $|P_n| = \infty$ for every $n \in \nn$. By our initial observation, we can construct a sequence $(M_n)_{n \in \nn}$ of non-finitely generated Puiseux monoids having full systems of sets of lengths such that $\mathsf{d}(M_i^\bullet) \cap \mathsf{d}(M_j^\bullet) = \{1\}$ when $i \neq j$. If $M_i \cong M_j$ for $i, j \in \nn$, then \cite[Proposition~3.2]{fG18} ensures that $M_j = rM_i$ for some $r \in \qq_{> 0}$. If $i \neq j$, then taking $q \in M_j^\bullet$ with $\mathsf{d}(q) > \mathsf{d}(r)$ we would obtain that $q \notin r M_i$. Hence $i = j$. Thus, $(M_n)_{n \in \nn}$ consists of pairwise non-isomorphic Puiseux monoids.
\end{proof}

\begin{prop}\emph{(cf.} \cite[Corollary~3.5]{GS18}\emph) \label{prop:infinitely many PA with full system of sets of lengts}
	There are infinitely many non-isomorphic Puiseux algebras having full systems of sets of lengths.
\end{prop}

\begin{proof}
	First, suppose that $M$ is a Puiseux monoid with full system of sets of lengths, which exists by Lemma~\ref{lem:PMs with full systems of sets of lengths}. Then $M$ is a BFM, and so $F[M]$ is a BFD by part~(1) of Proposition~\ref{prop:BF for PAlgebras}. Since the equality $\mathsf{L}_{F[M]}(X^q) = \mathsf{L}_M(q)$ holds for every $q \in M^\bullet$, it follows that $\mathcal{L}(M) \subseteq \mathcal{L}(F[M])$. Because $F[M]$ is a BFD, the fact that $M$ has full system of sets of lengths guarantees that $F[M]$ also has full system of sets of lengths. Take now an infinite class $\{M_i \mid i \in I\}$ of non-isomorphic Puiseux monoids with full systems of sets of lengths. Then the Isomorphism Problem for Puiseux algebras (Theorem 3.2) guarantees that $\{F[M_i] \mid i \in I\}$ is an infinite class of non-isomorphic Puiseux algebras having full systems of sets of lengths.
\end{proof}

The Characterization Problem for a class of atomic monoids $\mathcal{C}$ refers to the question of whether the function $M \mapsto \mathcal{L}(M)$ is injective on $\mathcal{C}$. Nontrivial instances of the Characterization Problem have been investigated in the past. For example, the Characterization Problem was answered negatively for the class of numerical monoids~\cite{ACHZ07}. Perhaps the most investigated instance of the Characterization Problem, which is still open, is for the class of Krull monoids with finite class group (Conjecture~\ref{conj:Characterization Problem for Krull monoids with finite class groups}). If $G$ is an additive finite abelian group with $|G| = n$, then the \emph{Davenport constant} of $G$, denoted by $D(G)$, is the minimum $d \in \nn$ such that every length-$d$ sequence of elements of $G$ contains a nonempty subsequence adding to zero.

\begin{conj} \label{conj:Characterization Problem for Krull monoids with finite class groups}
	Let $M$ and $M'$ be Krull monoids with respective finite abelian class groups $G$ and $G'$, each of their classes containing at least one prime divisor. Assume also that $D(G) \ge 4$. If $\mathcal{L}(M) = \mathcal{L}(M')$, then $M = M'$.
\end{conj}

As an immediate consequence of Proposition~\ref{prop:infinitely many PA with full system of sets of lengts}, one obtains that, for every field $F$, the answer of the Characterization Problem for the class of Puiseux algebras over $F$ is negative.

\begin{cor}
	For every field $F$, Puiseux algebras over $F$ are not determined by their systems of sets of lengths.
\end{cor}

\medskip
\subsection{The Finite Factorization Property}

As for satisfying the ACCP or the bounded factorization property, for an integral domain $R$ and a torsion-free monoid $M$ of type zero, Kim proved in~\cite[Theorem~3.25]{hK98} that $R[M]$ is an FFD if and only if $R$ is an FFD and $M$ is an FFM. However, a similar result when $M$ is reduced (instead of a type zero monoid) remains an open question.

We proceed to provide a class of Puiseux algebras that are FFDs but not UFDs. A Puiseux monoid is said to be \emph{increasing} if it can be generated by an increasing sequence of rationals. Increasing Puiseux monoids were first studied in~\cite{GG18} and have been recently considered in \cite{mBA20,BG20,fG19}.

\begin{prop} \label{prop:increasing PA are FFD}
	Let $F$ be a field, and let $M$ be an increasing Puiseux monoid. Then $F[M]$ is an FFD. In addition, $F[M]$ is a UFD if and only if $M \cong (\nn_0,+)$.
\end{prop}

\begin{proof}
	Because $0$ is not a limit point of $M$, it follows from \cite[Proposition~4.5]{fG19} that $M$ is a BFM. Therefore part~(1) of Proposition~\ref{prop:BF for PAlgebras} guarantees that $F[M]$ is a BFD. To verify that $F[M]$ is indeed an FFD, suppose towards a contradiction that there is an $f \in F[M] \setminus F$ such that $D_f := \{g \in F[M]^\ast  \mid \,  g \! \mid_{F[M]} f\}$ contains infinitely many non-associate divisors of $f$. Since $M$ is increasingly generated, the set $M \cap (0,\deg f]$ is finite. Clearly, for each $g \in D_f$, the inclusion $\supp(g) \subseteq M \cap (0,\deg f]$ holds. Hence there exists $S \subseteq M \cap (0,\deg f]$ such that the set $\{g \in D_f \mid \supp(g) = S\}$ contains infinitely many non-associate divisors of~$f$. Let $m$ be the least common multiple of $\mathsf{d}(M \cap (0, \deg f])$. Observe that
	\[
		G = \{g(X^m) \mid g \in D_f \ \text{and} \ \supp(g) = S \}
	\]
	is a subset of $F[X]$ consisting of infinitely many divisors of $f(X^m)$ in $F[X]$. Because $F[X]$ is a UFD, there exists $g_1$ and $g_2$ in $D_f$ such that the elements $g_1(X^m)$ and $g_2(X^m)$ of $G$ are associates in $F[X]$. As $F[M]^\times = F^\times = F[X]^\times$, it follows that $g_1$ and $g_2$ must be associates in $F[M]$, which is a contradiction. Thus, each element of $F[M]$ has only finitely many non-associate divisors and, because $F[M]$ is atomic, it is an FFD by \cite[Theorem~5.1]{AAZ90}. The fact that $F[M]$ is a UFD if and only if $M \cong (\nn_0,+)$ follows from \cite[Theorem~4.2]{fG20}.
\end{proof}

The Puiseux monoids $S_r := \langle r^n \mid n \in \nn_0 \rangle$, where $r \in \qq_{> 0}$, have been recently studied in~\cite{CGG20} under the term \emph{cyclic rational semirings} (clearly, they are closed under multiplication). It is known that for every $r \in \qq_{>0}$ with $r \notin \nn \cup \{1/n \mid n \in \nn\}$, the monoid $S_r$ is atomic with $\mathcal{A}(S_r) = \{r^n \mid n \in \nn_0\}$ (see \cite[Proposition~4.3]{CGG20a}).
 
\begin{cor} \label{cor:PA of cyclic rational semirings are FFD when r>1}
	 For each $r \in \qq_{\ge 1}$, the Puiseux algebra $F[S_r]$ is an FFD, and it is a UFD if and only if $r \in \nn$.
\end{cor}

\begin{proof}
	The fact that $F[S_r]$ is an FFD follows as a direct consequence of Proposition~\ref{prop:increasing PA are FFD} as $S_r$ is an increasing monoid when $r \in \qq_{\ge 1}$. For the second statement, it is clear that $S_r = \nn_0$ when $r \in \nn$. On the other hand, if $r = a/b$, where $a,b \in \nn_{\ge 2}$ and $\gcd(a,b) = 1$, then because $a = a \cdot 1 = b \cdot r$ the element $a$ has two distinct factorizations in $S_r$, and so $S_r \not\cong (\nn_0, +)$. Thus, $S_r \cong (\nn_0, +)$ if and only if $r \in \nn$, and so the last statement of the corollary follows from the last statement of Proposition~\ref{prop:increasing PA are FFD}.
\end{proof}

\smallskip
\subsection{Further Observations}

We have seen before that a Puiseux algebra satisfies the ACCP (resp., is a BFD) if and only if its exponent Puiseux monoid satisfies the ACCP (resp., is a BFM). For a general torsion-free monoid $M$, there seems to be no characterization (in terms of $M$) for the monoid algebra $F[M]$ to satisfy the ACCP, being a BFD, or being an FFD (see \cite[page 34]{AG20} for more details). In addition, it seems to be still open whether a monoid algebra $F[M]$ is an FFD provided that $M$ is a torsion-free reduced FFM.

In this section, we have constructed Puiseux algebras witnessing the failure of the reverse statements of all the implications in Diagram~\eqref{diag:AAZ chain}, except the last one. Although we still do not know whether the last implication is reversible, we suspect it is not. We finish this section proposing a Puiseux algebra as a potential counterexample. For $r \in (0,1) \cap \qq$ with $\mathsf{n}(r) \neq 1$, consider the Puiseux monoid $S_r$. Since
\[
	\mathsf{n}(r)r^n = \mathsf{d}(r)r^{n+1} = (\mathsf{d}(r) - \mathsf{n}(r))r^{n+1} + \mathsf{n}(r)r^{n+1}
\]
for every $n \in \nn$, the sequence $(\mathsf{n}(r)r^n + S_r)_{n \in \nn}$ is an ascending chain of principal ideals of~$S_r$. Clearly, this sequence does not stabilize, and so $S_r$ does not satisfy the ACCP. Now, let $F$ be a field. By \cite[Proposition~1.4]{hK01}, the Puiseux algebra $F[S_r]$ does not satisfy the ACCP. However, we believe that $F[S_r]$ is an atomic domain. The case when $r \ge 1$ in the following conjecture follows from Corollary~\ref{cor:PA of cyclic rational semirings are FFD when r>1}.

\begin{conj} \label{conj:an atomic PA}
	Let $F$ be a field, and take $r \in \qq_{>0}$. If $S_r$ is an atomic monoid, then $F[S_r]$ is an atomic domain. 
\end{conj}

\bigskip
\section{Classes of Antimatter Puiseux Algebras}
\label{sec:antimatter PA}

In this section we prove that the seminormal closure, root closure, and complete integral closure of a Puiseux algebra are equal, and we describe such closures in terms of the exponent Puiseux monoid. Our description will yield a class of antimatter and seminormal Puiseux algebras. We will also offer another class of antimatter Puiseux algebras that are not seminormal. Before proceeding, we would like to emphasize that antimatter domains were first investigated in~\cite{CDM99} and classes of antimatter Puiseux algebras were first constructed in~\cite{ACHZ07}.

\medskip
\subsection{Algebraic Closures}

Let $R$ be an integral domain with quotient field denoted by $\qf(R)$. The \emph{seminormal closure}, \emph{root closure}, and \emph{complete integral closure} of $R$, respectively denoted by $R'$, $\bar{R}$, and $\widehat{R}$, are the overrings of $R$ whose multiplicative monoids are ${R^*}'$, $\overline{R^*}$, and $\widehat{R^*}$, respectively. Thus,
\begin{equation} \label{eq:closure containments}
	R \subseteq R' \subseteq \bar{R} \subseteq \widehat{R} \subseteq \qf(R).
\end{equation}
The integral domain $R$ is called \emph{seminormal} (resp., \emph{root closed} or \emph{completely integrally closed}) if $R' = R$ (resp., $\bar{R} = R$ or $\widehat{R} = R$). In general, $R' \neq \bar{R}$ and $\bar{R} \neq \widehat{R}$ even in the context of monoid algebras.
 
 \begin{example} \hfill
 	\begin{enumerate}
 		\item In~\cite[Example~2.56]{BG09}, W. Bruns and J. Gubeladze exhibit an additive submonoid $M$ of $\nn_0^5$ that is seminormal but not root closed. Since $\qq[M]' = \qq[M']$ by \cite[Corollary~4.77]{BG09} and $\overline{\qq[M]} = \qq[\bar M]$ by \cite[Corollary~12.11]{rG84}, one obtains that $\qq[M]' \neq \overline{\qq[M]}$. 
 		\smallskip
 		
 		\item Consider the additive submonoid $M := \{ (0,0) \} \cup \nn^2$ of $\nn_0^2$, which satisfies that $\gp(M) = \zz^2$. It follows immediately that $M$ is root closed, and therefore, \cite[Corollary~12.11]{rG84} guarantees that the monoid algebra $\qq[M]$ is also root closed. Notice, on the other hand, that $M$ is not completely integrally closed because $(1,1) + n(0,1) \in M$ for every $n \in \nn$ even though $(0,1) \notin M$. So it follows from \cite[Corollary~12.7]{rG84} that $\qq[M]$ is not completely integrally closed. As a result, $\overline{\qq[M]} \neq \widehat{\qq[M]}$. 
 	\end{enumerate}
 \end{example}
 
 However, as we shall prove in the next theorem, in the class consisting of Puiseux algebras the three algebraic closures above coincide. First, let us argue the following lemma.
 
 \begin{lemma} \label{lem:intersection of monoid algebra and quotient field}
 	Let $F$ be a field, and let $M$ be a Puiseux monoid. Then the equality $F[M] \cap F(X) = F[X]$ holds.
 \end{lemma}

\begin{proof}
	It suffices to argue that $F[M] \cap F(X) \subseteq F[X]$, as the reverse inclusion follows immediately. To do this, take $f = \sum_{i=1}^k \alpha_i X^{q_i} \in F[M] \cap F(X)$ represented in canonical form as an element of $F[M]$. Let $\ell$ be the least common multiple of $\mathsf{d}(q_1), \dots, \mathsf{d}(q_k)$. Then take $g = \sum_{i=1}^m \beta_i X^{m_i}$ and $h = \sum_{i=1}^n \theta_i X^{n_i}$ in $F[X]$, both of them represented in canonical form, such that $f = g/h$. Then
	\begin{equation} \label{eq:F[M] cap F(X)}
		\sum_{i=1}^n \theta_i X^{\ell n_i} \sum_{i=1}^k \alpha_i X^{k_i} = h(X^\ell) f(X^\ell) = g(X^\ell) = \sum_{i=1}^m \beta_i X^{\ell m_i},
	\end{equation}
	in $F[X]$, where $k_i := \ell q_i \in \nn$ for every $i \in \ldb 1,k \rdb$. Let us argue inductively that $q_i \in \nn$ for every $i \in \ldb 1,k \rdb$. As $k_1 = \ell m_1 - \ell n_1$, we see that $q_1 \in \nn$. Suppose that $q_1, \dots, q_j \in \nn$ for some $j \in \ldb 1,k-1 \rdb$. Consider the monomial $\theta_1 \alpha_j X^{\ell n_1 + k_j}$ that shows when one multiplies out the leftmost part of \eqref{eq:F[M] cap F(X)}. If $\ell n_1 + k_j \in \supp \, g(X^\ell)$, then $\ell$ must divide $\ell n_1 + k_j$, and therefore, $q_j \in \nn$. If $\ell n_1 + k_j \notin \supp \, g(X^\ell)$, then the monomial $\theta_1 \alpha_j X^{\ell n_1 + k_j}$ should cancel with monomials of the form $\theta_i \alpha_t X^{\ell n_i + k_t}$ with $t < j$, in which case $\ell$ must divide $k_j - k_t$. As $\ell$ divides $k_t$, it follows that $q_j \in \nn$. Then we conclude that $f \in F[X]$. Thus, $F[M] \cap F(X) \subseteq F[X]$.
\end{proof}

\begin{theorem} \label{thm:closures of PA are antimatter over algebraically closed fields}
	Let $F$ be a field, and let $M$ be a Puiseux monoid. Then the following statement hold.
	\begin{enumerate}
		\item $F[M]' = \overline{F[M]} = \widehat{F[M]} = F[\emph{\gp}(M) \cap \qq_{\ge 0}]$.
		\smallskip
		
		\item If $M$ is finitely generated, then $F[M]'$ is atomic.
		\smallskip
		
		\item If $F$ is algebraically closed and $M$ is not finitely generated, then $F[M]'$ is antimatter.
	\end{enumerate}
\end{theorem}

\begin{proof}
	(1) By virtue of~\eqref{eq:closure containments}, it suffices to argue that $\widehat{F[M]} \subseteq F[\gp(M) \cap \qq_{\ge 0}]$ and $F[\gp(M) \cap \qq_{\ge 0}] \subseteq F[M]'$. To verify the latter inclusion, it is enough to observe that the equality $\gp(M) \cap \qq_{\ge 0} = M'$ holds by~\cite[Proposition~3.1]{GGT21} while the equality $F[M'] = F[M]'$ holds by \cite[Corollary 4.77]{BG09}.
	
	To prove the former inclusion, take $f$ in the complete integral closure of $F[M]$, and then take $g \in F[M]$ such that $gf^n \in F[M]$ for every $n \in \nn$. Write $f = f_1/f_2$ for $f_1, f_2 \in F[M]$ with $f_2 \neq 0$. Assume, by way of contradiction, that $f_2$ does not divide $f_1$ in $F[\gp(M) \cap \qq_{\ge 0}]$. Now let $\ell$ be the least common multiple of the set $\bigcup_{s \in S} \mathsf{d}(\supp(s))$, where $S := \{g, f_1, f_2\}$. It is clear that $n/\ell \in \gp(M) \cap \qq_{\ge 0}$ for every $n \in \nn$. Therefore setting $q(X) := f_1(X^\ell)/f_2(X^\ell) \in F(X)$, one can see that $q \notin F[X]$ as otherwise $f_1/f_2 = q(X^{1/\ell}) \in F[\gp(M) \cap \qq_{\ge 0}]$, which is not possible. Then $f_2(X^\ell)$ does not divide $f_1(X^\ell)$ in $F[X]$. As $F[X]$ is a UFD, there exist $a \in \ii(F[X])$ and $m \in \nn$ such that
	\[
		a(X)^m \mid_{F[X]} f_2(X^\ell) \quad \text{but} \quad a(X)^m \nmid_{F[X]} f_1(X^\ell).
	\]
	Let $\mu = \max\{n \in \nn \mid a(X)^n \mid_{F[X]} f_1(X^\ell)\}$. Take $N \in \nn$ such that the inequality $N \deg a > \ell \deg g$ holds, and then take $h_N \in F[M]$ such that $g f_1^N = f_2^N h_N$. Observe that $h_N(X^\ell) = g(X^\ell)f_1(X^\ell)^N f_2(X^\ell)^{-N} \in F(X)$. Then it follows from Lemma~\ref{lem:intersection of monoid algebra and quotient field} that $h_N(X^\ell) \in F[X]$. As a result, the factors in $g(X^\ell) f_1(X^\ell)^N = f_2(X^\ell)^N h_N(X^\ell)$ are polynomials in $F[X]$. This, together with the fact that $a(X)^m \mid_{F[X]} f_2(X^\ell)$, implies that $a(X)^{mN} \mid_{F[X]} g(X^\ell) f_1(X^\ell)^N$. So there exists $b(X) \in F[X]$ such that
	\[
		b(X)a(X)^{N(m-\mu)} = g(X^\ell) \bigg( \frac{f_1(X)^\ell}{a(X)^\mu} \bigg)^N.
	\]
	Since $m > \mu$ and $F[X]$ is a UFD, $a(X)^N \mid_{F[X]} g(X^\ell)$. However, this contradicts the inequality $N \deg a > \ell \deg g$. Therefore $f \in F[\gp(M) \cap \qq]$.  We conclude that $\widehat{F[M]} \subseteq F[\gp(M) \cap \qq_{\ge 0}]$.
	\smallskip
	
	(2) Suppose that $M$ is finitely generated, namely, $M = \langle q_1, \dots, q_n \rangle$ for some $n \in \nn$ and $q_1, \dots, q_n \in \qq_{> 0}$. Letting $\ell$ be the least common multiple of $\mathsf{d}(q_1), \dots, \mathsf{d}(q_n)$ and~$g$ be the greatest common divisor of $\mathsf{n}(q_1), \dots, \mathsf{n}(q_n)$, one can check that $N := \ell g^{-1}M$ is a numerical monoid. Therefore $N' = \nn_0$. It is clear that $M \cong N$. Then it follows from \cite[Corollary 4.77]{BG09} and Theorem~\ref{thm:the isomorphism problem} that $F[M]' = F[M'] \cong F[N'] \cong F[X]$, and so $F[M]'$ is a UFD and, in particular, an atomic domain.
	\smallskip
	
	(3) Suppose now that $M$ is not finitely generated. In light of Theorem~\ref{thm:the isomorphism problem}, one can replace $M$ by $(\gcd \mathsf{n}(M^\bullet))^{-1}M$ and assume that $\gcd \mathsf{n}(M^\bullet) = 1$. Then \cite[Proposition~3.1]{GGT21} ensures that $M' = \langle \mathsf{d}(q)^{-1} \mid q \in M^\bullet \rangle$. From the fact that $M$ is not finitely generated, one can deduce that $|\mathsf{d}(M'^\bullet)| = |\mathsf{d}(M^\bullet)| = \infty$. We check that $M'$ is pure (i.e., for each $b \in M'$ there exists $n \in \nn_{\ge 2}$ such that $b/n \in M'$). To do so, take $b \in M'^\bullet$. As $|\mathsf{d}(M^\bullet)| = \infty$, we can take $d \in \mathsf{d}(M^\bullet)$ such that $d \nmid \mathsf{d}(b)$. It is clear that the least common multiple $\ell$ of $\mathsf{d}(b)$ and $d$ belongs to $\mathsf{d}(M'^\bullet)$. Setting $n := \ell/\mathsf{d}(b)$, we obtain that $n \ge 2$ and $b/n = \mathsf{n}(b)/\ell \in M'$. Hence $M'$ is pure. Since $F$ is algebraically closed, it follows from \cite[Corollary 4.77]{BG09} and \cite[Theorem~1]{ACHZ07} that $F[M]' = F[M']$ is an antimatter domain.
\end{proof}

Recall that for each $r \in \qq_{>0}$, the Puiseux monoid $\langle r^n \mid n \in \nn_0 \rangle$ is denoted by $S_r$.

\begin{cor} \label{cor:infinitely many PA over algebraically closed fields}
	Let $F$ be an algebraically closed field. For each $p \in \pp$, the Puiseux algebra $F[S_{1/p}]$ is antimatter. In addition, $F[S_{1/p}] \not\cong F[S_{1/q}]$ if $q \in \pp \setminus \{p\}$.
\end{cor}

\begin{proof}
	As $-S_{1/p} \cup S_{1/p}$ is an additive subgroup of $\qq$, the equality $S_{1/p} = \gp(S_{1/p}) \cap \qq_{\ge 0}$ holds, and so $F[S_{1/p}] = F[\gp(S_{1/p}) \cap \qq_{\ge 0}] = F[S_{1/p}]'$ by part~(1) of Theorem~\ref{thm:closures of PA are antimatter over algebraically closed fields}. Since~$M$ is not finitely generated and $F$ is algebraically closed, it follows from part~(3) of Theorem~\ref{thm:closures of PA are antimatter over algebraically closed fields} that $F[S_{1/p}]$ is antimatter.
	
	To argue the second statement, suppose for the sake of a contradiction that there exist $p,q \in \pp$ with $p \neq q$ such that $F[S_{1/p}] \cong F[S_{1/q}]$. It follows from Theorem~\ref{thm:the isomorphism problem} that $S_{1/p} \cong S_{1/q}$, and therefore, \cite[Proposition~3.2]{fG18} guarantees that $S_{1/q} = r S_{1/p}$ for some $r \in \qq_{> 0}$. Taking $n \in \nn$ such that $p^n \nmid \mathsf{n}(r)$, one obtains that $r/p^n \in rS_{1/p} = S_{1/q}$, which contradicts that $\mathsf{d}(S_{1/q}) = \{q^n \mid n \in \nn_0\}$.
\end{proof}
\smallskip

With the notation as in Corollary~\ref{cor:infinitely many PA over algebraically closed fields}, the hypothesis that the field $F$ is algebraically closed is not superfluous, as we will confirm in Example~\ref{ex:antimatter and idf are not transfer properties}.

A monoid is an \emph{irreducible-divisor-finite monoid} (or an \emph{IDFM}) if each element is divisible by only finitely many atoms up to associates, while we say that an integral domain is an \emph{irreducible-divisor-finite domain} (or an \emph{IDFD}) if its multiplicative monoid is an IDFM. It has been proved in \cite[Theorem~5.1]{AAZ90} that an FFD can be characterized by being an atomic IDFD; this result was generalized for monoids in~\cite[Theorem~2]{fHK92}. Although we do not know whether a Puiseux algebra $F[M]$ is an FFD provided that $M$ is an FFM, we can answer the corresponding question for the irreducible-divisor-finite property. In the following example we verify that a Puiseux algebra $F[M]$ need not be antimatter (resp., an IDFD) when~$M$ is antimatter (resp., an IDFM).

We let $\Phi_n(X)$ denote the $n$-th cylcotomic polynomial, while we let $\varphi$ denote Euler's totient function. 

\begin{example} \label{ex:antimatter and idf are not transfer properties}
	Take $p \in \pp$, and consider the Puiseux monoid $S_{1/p}$. We have seen that $S_{1/p}$ is antimatter and, therefore, an IDFM. To argue that the Puiseux algebra $\qq[S_{1/p}]$ is neither antimatter nor an IDFD, it suffices to show that the element $X - 1 \in \qq[S_{1/p}]$ is divisible by infinitely many non-associate irreducible elements in~$\qq[S_{1/p}]$. For every $k \in \nn$, we can factor $X-1$ as follows:
	\begin{equation} \label{eq:factorization into cyclotomic polynomials}
		X - 1 = \big( X^{\frac{1}{p^k}} \big)^{p^k}  - 1 = \prod_{n=0}^k \Phi_{p^n} \big( X^{\frac{1}{p^k}} \big).
	\end{equation}
	We claim that each factor on the rightmost expression of \eqref{eq:factorization into cyclotomic polynomials} is irreducible. To verify that our claim holds, fix $n \in \ldb 0,k \rdb$ and write $\Phi_{p^n}\big( X^{1/{p^k}} \big) = f(X)g(X)$, where $f(X), g(X) \in \qq[S_{1/p}]$. Let $p^m$ be the least common multiple of $\mathsf{d} \big( \supp(f) \cup \supp(g) \big)$. Since $\deg \Phi_{p^n}\big( X^{1/{p^k}}\big) = \varphi(p^n)/p^k = \frac{p-1}{p^{k-n+1}}$, it follows that $p^{k-n+1}$ divides $p^m$ and, therefore, $ m-k+n \ge 1$. Because
	\[
		\Phi_{p^n}\big( X^{p^{m-k}} \big) = \Phi_p \big( X^{p^{m-k+n-1}} \big) = \Phi_{p^{m-k+n}} (X),
	\]
	it follows that $\Phi_{p^n}\big( X^{p^{m-k}} \big)$ is an irreducible polynomial in $\qq[X]$. Now the fact that $\Phi_{p^n}\big( X^{p^{m-k}} \big) = f(X^{p^m})g(X^{p^m})$ yields a factorization of $\Phi_{p^n}\big( X^{p^{m-k}} \big)$ in $\qq[X]$ implies that either $f(X^{p^m}) \in \qq$ or $g(X^{p^m}) \in \qq$, which in turns implies that either $f(X) \in \qq$ or $g(X) \in \qq$. So $\Phi_{p^n}\big( X^{1/{p^k}} \big)$ is irreducible in the Puiseux algebra $\qq[S_{1/p}]$. Because $k$ was taking arbitrarily in $\nn$ and $\Phi_{p^{n_1}}\big( X^{1/{p^k}} \big) \neq \Phi_{p^{n_2}}\big( X^{1/{p^k}} \big)$ whenever $n_1 \neq n_2$, we conclude that $X-1$ has infinitely many irreducible divisors in $\qq[S_{1/p}]$. Hence $\qq[S_{1/p}]$ is neither antimatter nor an IDFD.
\end{example}

In the direction of Corollary~\ref{cor:infinitely many PA over algebraically closed fields}, we have the following question.

\begin{question}
	Is there an antimatter Puiseux monoid that is not root closed such that the algebra $F[M]$ is antimatter over any (or some) algebraically closed field $F$?
\end{question}
\smallskip

The antimatter Puiseux algebras we have seen so far come from part~(3) of Theorem~\ref{thm:closures of PA are antimatter over algebraically closed fields} and are, therefore, seminormal. By contrast, we would like to construct a class of antimatter Puiseux algebras that are not seminormal. For distinct $p,q \in \pp$, let $M_{p,q}$ denote the Puiseux monoid $\langle p^{-m}q^{-n} \mid m,n \in \nn_0 \rangle$.

\begin{prop} \label{prop:characteristic p and p-divisible implies that F[M] is antimatter}
	Let $F$ be a perfect field of finite characteristic $p$. For each $q \in \pp \setminus \{p\}$, the Puiseux algebra $F[M_{p,q}]$ is antimatter but fails to be seminormal. In addition, $F[M_{p,q}] \not\cong F[M_{p,q'}]$ for any $q' \in \pp \setminus \{p,q\}$.
\end{prop}

\begin{proof}
	Fix $q \in \pp \setminus \{p\}$. For each $x \in M_{p,q}$ it is clear that $x/p \in M_{p,q}$, and therefore, $M_{p,q}$ is antimatter. To argue that $F[M_{p,q}]$ is an antimatter domain, consider the element
	\[
		f = \alpha_1 X^{q_1} + \dots + \alpha_n X^{q_n} \in F[M_{p,q}] \setminus F.
	\]
	As $F$ is a perfect field of characteristic $p$, the Frobenius homomorphism $x \mapsto x^p$ is surjective, and so for each $i \in \ldb 1,n \rdb$, there exists $\beta_i \in F$ with $\alpha_i = \beta_i^p$ for some $\beta_i \in F$. In addition, it follows from our initial observation that $q_i/p \in M_{p,q}$ for every $i \in \ldb 1,n \rdb$. Therefore the element $f = \big( \beta_1 X^{q_1/p} + \dots + \beta_n X^{q_n/p} \big)^p$ is not irreducible in $F[M_{p,q}]$. Hence $F[M_{p,q}]$ is an antimatter Puiseux algebra.
	
	In light of \cite[Corollary~4.77]{BG09}, proving that the Puiseux algebra $F[M_{p,q}]$ is not seminormal amounts to verifying that $M_{p,q}$ is not a seminormal monoid. Assume, by way of contradiction, that $M$ is seminormal. Then $\frac{1}{pq} \in M$ by \cite[Proposition~3.1]{GGT21}. So we can write
	\begin{align} \label{eq:antimatter not ascending cyclic}
		\frac{1}{pq} = \sum_{i=1}^t \frac{\alpha_i}{p^i} + \sum_{i=1}^s \frac{\beta_i}{q^i}
	\end{align}
	for some coefficients $\alpha_1, \dots \alpha_t, \beta_1, \dots, \beta_s \in \nn_0$ with $\alpha_t \neq 0$ and $\beta_s \neq 0$. After simplifying if necessary, we can assume that $p \nmid \alpha_i$ and $q \nmid \beta_j$ for any $i \in \ldb 1,t \rdb$ and $j \in \ldb 1,s \rdb$. Multiplying (\ref{eq:antimatter not ascending cyclic}) by $p^t q^s$ one obtains that $t = s = 1$. However, $\frac{\alpha_1}{p} + \frac{\beta_1}{q} \ge \frac{1}{p} + \frac 1q > \frac{1}{pq}$, which contradicts~\eqref{eq:antimatter not ascending cyclic}. Hence $F[M_{p,q}]$ is not a seminormal domain.
	
	To argue that $F[M_{p,q}] \not\cong F[M_{p,q'}]$ for any $q' \in \pp \setminus \{p,q\}$ one can merely mimic the lines of the second paragraph of the proof of Corollary~\ref{cor:infinitely many PA over algebraically closed fields}.
\end{proof}

\begin{remark} 
	Proposition~\ref{prop:characteristic p and p-divisible implies that F[M] is antimatter} is a version of \cite[Theorem~5(2)]{ACHZ07}, which states that if~$R$ is an antimatter GCD-domain whose quotient field is perfect of finite characteristic, then $R[\qq_{\ge 0}]$ is also an antimatter GCD-domain.
\end{remark}

\bigskip
\section*{Acknowledgments}
	The author would like to thank an anonymous referee for useful suggestions that help to improve the final version of this paper. During the preparation of the same, the author was supported by an NSF-AGEP graduate fellowship and by the NSF postdoctoral award DMS-1903069.

\bigskip

\bigskip

\end{document}